\documentclass[reqno]{amsart}
\usepackage{graphicx}
\usepackage{epsfig}
\usepackage{amssymb}
\usepackage{xcolor}
\usepackage[
backend=biber,
style=alphabetic,
sorting=nyt
]{biblatex}
\usepackage{hyperref}
\usepackage{amsfonts}
\usepackage{dsfont}
\usepackage{amsmath}
\usepackage{stmaryrd}
\usepackage{subfigure, verbatim}
\newtheorem{theorem}{Theorem}
\newtheorem*{theorem*}{Theorem}
\newtheorem{lemma}{Lemma}
\newtheorem{proposition}{Proposition}[section]
\theoremstyle{definition}
\newtheorem{definition}{Key Proposition}

\newtheorem{corollary}{Corollary}

\theoremstyle{remark}
  
\numberwithin{equation}{section}
\newcommand\numberthis{\addtocounter{equation}{1}\tag{\theequation}}
\theoremstyle{remark}

\newcommand{\beq}{\begin{equation}}
\newcommand{\eeq}{\end{equation}}
\newcommand{\beqn}{\begin{equation*}}
\newcommand{\eeqn}{\end{equation*}}
\newcommand{\bea}{\begin{eqnarray}}
\newcommand{\eea}{\end{eqnarray}}
\newcommand{\bean}{\begin{eqnarray*}}
\newcommand{\eean}{\end{eqnarray*}}

\newcommand{\be}{\begin{enumerate}}
\newcommand{\ee}{\end{enumerate}}
\newcommand{\bi}{\begin{itemize}}
\newcommand{\ei}{\end{itemize}}
\newcommand{\bd}{\begin{description}}
\newcommand{\ed}{\end{description}}

\title[Almost sure upper bounds]{Two almost sure upper bounds for sums of random multiplicative functions}
\author{Christopher Atherfold}
\address{School of Mathematics, University of Bristol, Bristol, BS8 1UG, United Kingdom}
\email{christopher.atherfold@bristol.ac.uk}
\begin{document}
\begin{abstract}
    Let $f,g$ be a squarefree and extended Rademacher random multiplicative functions, respectively. We prove that, for every $\varepsilon>0$, almost surely
    $$\big|\sum_{n \le x}\frac{f(n)}{\sqrt{n}}\big| \ll (\log\log x)^{1/4+\varepsilon}, \qquad \big|\sum_{n\le x}g(n)\big|\ll\sqrt{x}(\log\log x)^{1/4+\varepsilon}.$$
    We also show that almost surely, there exist arbitrarily large values of $x$ such that $\sum_{n \le x}\frac{f(n)}{\sqrt{n}} \gg (\log\log x)^{1/4+o(1)}$. We follow Caich and Durkan and Pearce-Crump by using a version of the Azuma--Hoeffding inequality for martingale difference sequences.
\end{abstract}
\maketitle
\section{Introduction}

\subsection{Motivation and Results}
Let $f(n)$ be a Rademacher random multiplicative function, so we consider the sequence $(f(p))$ to be a sequence of independent identically distributed Rademacher random variables on $\{\pm 1\}$ and extend these multiplicatively to all squarefree natural numbers, and set $f(n) =0$ otherwise. We also define $g(n)$ to be an extended Rademacher random multiplicative function if instead the $f(p)$ are extended completely multiplicatively to all natural numbers instead (note $f(n)=g(n)$ when $n$ is squarefree). We denote $\mathbb{E}$ as the expectation taken over all possible realisations of $f$ (which are determined purely by the values of $f(p)$). \par
In this paper we are concerned with upper bounds of the large fluctuations of the partial sums
$$\sum_{n \le x}\frac{f(n)}{\sqrt{n}}, \qquad \sum_{n \le x}g(n).$$
The Rademacher random multiplicative function was proposed by Wintner \cite{wintner1944} as a random model for the M\"obius function, with further developments due to Hal\'asz in \cite{halasz83} and Lau--Tenenbaum--Wu \cite{LTW13}. The recent developments in this area have been spearheaded by Harper, in particular the striking resolution of Helson's conjecture in \cite{aharper201}, which showed that
$$\mathbb{E}|\sum_{n \le x}f(n)| \asymp \frac{\sqrt{x}}{(\log\log x)^{1/4}}$$
through a connection to critical Gaussian multiplicative chaos. In particular, such sums have very heavy tails and are not distributed in a Gaussian sense. Despite these heavy tails, Harper proved in \cite{harper23} that for function $V(x) \to \infty$ that there almost surely exist arbitrarily large values of $x$ such that
$$|\sum_{n \le x}f(n)| \ge \frac{\sqrt{x}(\log\log x)^{1/4}}{V(x)},$$
showing that for large fluctuations, the Gaussian heuristic was not too far from the truth. There has been significant progress made on the upper bounds of the large fluctuations as well recently. Caich proved in \cite{caich23} that for any $\varepsilon>0$, then almost surely
$$\sum_{n \le x}f(n) \ll \sqrt{x}(\log\log x)^{3/4+\varepsilon}.$$
Following Caich's method, Durkan and Pearce-Crump refined this to $1/4+\varepsilon$ in \cite{dpc26}, and the state of the art is due to Verreault in \cite{Verreault26}, where he shows that for any $\varepsilon>0$, then almost surely
$$\sum_{n \le x}f(n) \ll \sqrt{x}(\log\log x)^{1/4}(\log\log\log x)^{1+\varepsilon}.$$
Verreault's method is adapted far more closely to the multiplicative chaos object connected to the sum, whereas Caich's method only requires a sufficiently strong bound on the low fractional moments of integrals involving random Euler products. \par
Weighted Rademacher and extended Rademacher expressions have been significantly less explored until recently in their moment and large fluctuation statistics. The author with Gerspach and Hamdan proved in \cite{AGH26} that one has both
$$\mathbb{E}|\sum_{n \le x}g(n)|^{2q} \asymp (\frac{x}{\sqrt{\log\log x}})^q, \qquad \mathbb{E}|\sum_{n \le x}g(n)| \asymp \sqrt{x}(\log\log x)^{3/4}$$
where $q \in (0,1/2)$ fixed (with similar results obtained in the weighted case). This is rather surprising given that from the contributions of the square numbers alone, one has expectation $\mathbb{E}[\sum_{n \le x}g(n)] = \lfloor\sqrt{x}\rfloor$. We note that the weighted Steinhaus sums were studied in the work of Gerspach for moments \cite{gerspach20} and Hardy for large fluctuations \cite{shardy23} respectively. Separately, Gou--Xu \cite{GX26} (and independently Hardy--Klurman \cite{HK}) proved that for any $V(x) \to \infty$, there exist arbitrarily large values of $x$ such that
$$\sum_{n \le x}g(n) \ge \frac{\sqrt{x}(\log\log x)^{1/4}}{V(x)}.$$
In this paper, we show that the exponent of $1/4$ obtained in \cite{GX26} and \cite{HK} is sharp.
\begin{theorem}\label{thm:extended_upper}
    Let $g$ be an extended Rademacher random multiplicative function. Then for any $\varepsilon>0$, almost surely
    $$\sum_{n \le x}g(n) \ll \sqrt{x}(\log\log x)^{1/4+\varepsilon}.$$
\end{theorem}
By a very similar method, one can prove the following theorem.
\begin{theorem}\label{thm:weighted_upper}
    Let $f$ be a squarefree Rademacher random multiplicative function. Then for any $\varepsilon>0$, almost surely
    $$\sum_{n \le x}\frac{f(n)}{\sqrt{n}} \ll (\log\log x)^{1/4+\varepsilon}.$$
\end{theorem}
We also have an immediate corollary of Theorem 1 in \cite{harper23} for existence of large values in the weighted Rademacher sum.
\begin{theorem}\label{thm:weighted_existence}
    Let $f$ be a squarefree Rademacher random multiplicative function. Then for any $V(x) \to \infty$, there almost surely exist arbitrarily large $x$ such that
    $$\Big|\sum_{n \le x}\frac{f(n)}{\sqrt{n}}\Big| \ge \frac{(\log\log x)^{1/4}}{V(x)}.$$
\end{theorem}
\begin{proof}
    This is almost immediate. Indeed, we have
    $$\frac{1}{\sqrt{N}}\sum_{n \le N}f(n) = \sum_{n \le N}\frac{f(n)}{\sqrt{n}} -\sum_{n \le N-1}\frac{\sqrt{n+1}-\sqrt{n}}{\sqrt{N}}\Big(\sum_{m \le n }\frac{f(m)}{\sqrt{m}}\Big).$$
    Noting that $\sqrt{n+1}-\sqrt{n}$ is positive and provides a telescoping sum, we can take absolute values on both sides, and apply the triangle inequality and take the maximum value of the weighted sum. Consequently, we have
    $$|\frac{1}{\sqrt{N}}\sum_{n \le N}f(n)| \le \max_{n \le N}\Big|\sum_{m \le n}\frac{f(m)}{\sqrt{m}}\Big|(1+\sum_{n=1}^{N-1}\frac{\sqrt{n+1}-\sqrt{n}}{\sqrt{N}}) \le 2\max_{n \le N}\Big|\sum_{m \le n}\frac{f(m)}{\sqrt{m}}\Big|.$$
    Applying Harper's Theorem 1 in \cite{harper23} with $U(x) = V(x)/2$ yields the result.
\end{proof}
Using our upper bounds, we obtain the following using Theorem 1 in \cite{GX26} and Theorem \ref{thm:weighted_existence}.
\begin{corollary}
    \textit{Suppose $g$ is an extended Rademacher random multiplicative function. Then almost surely}
    $$\limsup_{x \to \infty}\frac{\log\Big(1+\frac{\big|\sum_{n \le x}g(n)\big|}{\sqrt{x}}\Big)}{\log \log \log x} = \frac{1}{4}.$$
    \textit{Similarly, if $f$ is a squarefree Rademacher random multiplicative function, then almost surely}
    $$\limsup_{x \to \infty}\frac{\log\Big(1+\big|\sum_{n \le x}\frac{f(n)}{\sqrt{n}}\big|\Big)}{\log \log \log x} = \frac{1}{4}.$$
\end{corollary}
\begin{proof}
    The upper bound is provided by Theorem \ref{thm:extended_upper}. As for the lower bound, choosing $V(x) = (\log\log x)^{\delta}$ for $\delta \in (0,1/4)$ in \cite{GX26} Theorem 1.1 shows that almost surely, there exist arbitrarily large $x$ such that
    $$\frac{\big|\sum_{n \le x}g(n)\big|}{\sqrt{x}} \ge (\log\log x)^{1/4-\delta}.$$
    Taking $\delta$ through the countable sequence $2^{-n-2}$ yields the corollary. An identical argument follows in the weighted Rademacher case, replacing the upper and lower bound inputs accordingly.
\end{proof}
Our second corollary follows from combining Theorems \ref{thm:extended_upper} and \ref{thm:weighted_upper} to obtain another large fluctuation result. To see this, we have the following identity
$$\sum_{n \le x}g(n) = \sum_{n \le x}f(n)\Big\lfloor\sqrt{x/n}\Big\rfloor=\sqrt{x}\sum_{n \le x}\frac{f(n)}{\sqrt{n}} - \sum_{n \le x}f(n)\Big\{\sqrt{x/n}\Big\}.$$
By the triangle inequality, we obtain the result below.
\begin{corollary}
    \textit{Let $\varepsilon >0$. Then, almost surely, for all large $x$,}
    $$\sum_{n \le x}f(n)\Big\{\sqrt{x/n}\Big\} \ll \sqrt{x}(\log\log x)^{1/4+\varepsilon}.$$
\end{corollary}
We then have that almost surely for random variable $S_x = \tfrac{|\sum_{n \le x}g(n)|}{\mathbb{E}|\sum_{n \le x}g(n)|}$, then we have
$$S_x \underset{x \to \infty}{\longrightarrow} 0 \quad a.s, \qquad \mathbb{E}[S_x] =1.$$
Consequently $S_x$ fails to be uniformly integrable and strengthens  the convergence of $S_x \to 0$ from convergence in probability observed in \cite{AGH26} to an almost sure convergence. This is rather different to that observed in the Steinhaus or unweighted squarefree Rademacher cases. In particular, for a given $g$, it almost surely stays below the order of magnitude given by the first moment's order of magnitude. \par
We end by discussing the weighted sum of extended Rademacher random multiplicative functions
$$\sum_{n \le x}\frac{g(n)}{\sqrt{n}} = \sum_{b \le \sqrt{x}}\frac{1}{b}\sum_{n \le x/b^2}\frac{f(n)}{\sqrt{n}}.$$
From Theorem \ref{thm:weighted_upper}, it immediately follows that for any $\varepsilon>0$, then almost surely
$$\sum_{n \le x}\frac{g(n)}{\sqrt{n}} \ll \log x(\log\log x)^{1/4+\varepsilon}.$$
However, in \cite{AGH26}, the authors show for $0 \le q \le 1/2$ fixed
$$\mathbb{E}\Big|\sum_{n \le x}\frac{g(n)}{\sqrt{n}}\Big|^{2q} \ll (\log x)^{q+2q^2},$$
which suggests perhaps the large fluctuations should be of order
$$\sqrt{\log x} \exp(\sqrt{2\log_2 x \log_4 x})$$
instead following a log normal heuristic, which is similar to the weighted Steinhaus case explored in \cite{shardy23}. Note that $\log_k$ denotes the $k$-fold iterated logarithm in this context. \par
The log normal heuristic follows from relating the weighted sum to integrals of random Euler products, and suggesting that the mass of such an integral should be coming from
$$|F_y(1/2)|^2 = \prod_{p \le y}\Big|1+\frac{f(p)}{p^{1/2}}\Big|^{2}$$
in the weighted extended Rademacher case (determined in the almost sure sense). If we expand this Euler product, 
$$\exp(-2\sum_{p \le x}\log(1-f(p)p^{-1/2})) = \exp\Big(\sum_{p \le x}(2f(p)p^{-1/2}+1/p + O(p^{-3/2}))\Big).$$
The $1/p$ factor introduces a deterministic drift, on top of the random walk heuristic for the random sum. In the weighted squarefree Rademacher case, then instead, one has a drift of $-1/p$. The negative drift makes the random Euler product at the central point significantly smaller, leading to the multiplicative chaos heuristic to prevail in the squarefree case. For more details on the weighted Rademacher case and its connections to \textit{multiplicative chaos with insertions}, see \cite{AGH26}. 
\subsection{A sketch of the proof}
Unlike when handling the low moments in \cite{AGH26}, the squares (or similar deterministic contributions) are not contributing to the term which is responsible for the large fluctuations. In particular, we have
$$\sum_{n \le x}g(n) = \lfloor\sqrt{x}\rfloor + \sum_{p \le x}f(p)\sum_{\substack{n \le x/p \\ P(\kappa(n)) < p}}g(n), \qquad \sum_{n \le x}\frac{f(n)}{\sqrt{n}} = 1+\sum_{p \le x}\frac{f(p)}{\sqrt{p}}\sum_{\substack{n \le x/p \\ P(n) < p}}\frac{f(n)}{\sqrt{n}}$$
where $\kappa(n)$ denotes squarefree part in the decomposition $n=\kappa(n)m^2$. From the definition of an extended Rademacher random multiplicative function, we rewrite the term on the right as
$$\sum_{p \le x}f(p)\sum_{\substack{n \le x/p \\ P(\kappa(n)) < p}}g(n) = \sum_{p \le x} A_p(x).$$
This is a martingale difference sequence with respect to the filtration $\mathcal{F}_{p^-} = \sigma\{f(q):q <p\}$. This is of the correct form to use a version of the Azuma--Hoeffding inequality adapted for martingale difference sequences (which from this point we will refer to as the Caich--Hoeffding inequality). From this point, the main point of the argument is to relate these partial sums to ones involving integrals of random Euler products and analyse how large these can get, so it becomes very important to understand the size of quantities like
$$\int_{-1/2}^{1/2}\left|\frac{F_y(1/2+ih)}{1/\log y +ih}\right|^2dh$$
with high probability, where we define for $f(p)$ independent Rademacher random variables
$$F_y(s) = \prod_{p \le y}(1+f(p)/p^{-s}).$$
It suffices to understand the low moments of such objects, which were recently established in \cite{AGH26}. This with some extra technical estimates is sufficient to obtain the bounds in the main theorems.
The difference between understanding weighted sums of squarefree Rademacher random multiplicative functions and unweighted extended Rademacher sums comes from a squarefree Rademacher sum with fractional part coefficients, so we go through the argument with the extended Rademacher case in mind, with suitable changes mentioned when handling the weighted squarefree case. \par
One thing which is rather nice in this proof is that we do not need to analyse any of the multiplicative chaos structures in too much detail, which considering the extended Rademacher case in particular is rather helpful. However, for sharper estimates it is likely one would need an approach similar to that in \cite{Verreault26}.
\subsection{Notation and organisation}
We need to set up and recap some notation here. Since we order the integers by their largest prime factors, we define $P(n)$ to be the largest prime factor of $n$, and let $P(1) = 1$ (so it is always included in the smooth sum). Further, recall that any integer can be written in the form $n = \kappa(n)m^2$ where $\kappa(n)$ is squarefree, and we use this distinction multiple times so we obtain the martingale difference sequence property in the extended Rademacher sum. Since we are considering martingales throughout, we let $\mathcal{F}_{p^-} = \sigma\{f(q) : q < p\}$ as defined in the previous section. We collect the auxiliary results in Section \ref{sec:prelim}, establish the various test point quantities we are investigating and discard one exceptional set in Section \ref{sec:setup} and the majority of the work occurs in Section \ref{sec:conditional_variance}.
\subsection*{Acknowledgements}
The author thanks Benjamin Durkan, Maxim Gerspach, Jad Hamdan, Seth Hardy and Besfort Shala for their thoughts, comments and suggestions relating to this project. He especially thanks the anonymous referee for their comments which significantly improved this paper. He also thanks his supervisors Oleksiy Klurman and Joseph Najnudel for their comments and encouragement on this project. The author is supported by the Heilbronn Institute of Mathematical Research.
\subsection*{AI usage declaration}
The author acknowledges the use of ChatGPT 6 Astra for discussions and proofreading surrounding this paper. All mathematical arguments developed with this assistance have been reworked and validated by hand and the author takes full responsibility for the accuracy and rigour of the proofs.
\section{Preliminary results}\label{sec:prelim}
\subsection{Number theoretic estimates}
For this upper bound, we will use much of the notation and framework established in Caich's proof of the unweighted Rademacher case. The first of these is the rough Rademacher hypercontractive inequality. This was first proved by Bonami \cite{bonami70} in a far more general setting than we present here, and reproved by Hal\'asz \cite{halasz83}.
\begin{proposition}
    (Probability Result 2.3, \cite{harper202}) Let $k \in \mathbb{N}$ and $f$ be a Rademacher multiplicative function. Then for $(a_n) \subset \mathbb{C}$ where $a_n \neq 0$ for only finitely many $n$, then
    \begin{equation*}
        \mathbb{E}\big|\sum_{n}a_nf(n)\big|^{2k} \le \big(\sum_n |a_n|^2\tau_{2k-1}(n)\big)^k 
    \end{equation*}
    where $\tau_{2k-1}(n) = \#\{(m_1,...,m_{2k-1}): m_1...m_{2k-1} = n\}$ is the $(2k-1)$-fold iterated divisor function. \label{prop_hyper}
\end{proposition}
We then need an estimate when transferring these statements to extended Rademacher random multiplicative functions on short intervals. We need the following estimates on divisor sums.
\begin{proposition}\label{prop:divisor_sums}
    (Tenenbaum (Ch. III, \cite{Tenenbaum15}; Shiu, \cite{Shiu80}) Let $K>0$. Then for any positive integer $k \le K$, then
    $$\sum_{n \le x}\tau_k(n) \ll x(\log x)^{k-1}, \qquad \sum_{n \le x}\frac{\tau_k(n)}{n} \ll (\log (2x))^k.$$
    Further, suppose $\sqrt{x} < h \le x$. Then
    $$\sum_{x < n \le x+h}\tau_{k}(n) \ll h(\log x)^{k-1}.$$
\end{proposition}
\begin{proposition}\label{prop:short_sum_count}
    Let $k \in \mathbb{N}$, $x > y$ and $g(n)$ an extended Rademacher random multiplicative function. Then
    $$\mathbb{E}\big|\sum_{x-y < n \le x}g(n)\big|^{2k} \ll \Big(\frac{y^2\log^{2k-1}(2x)}{x}+\sum_{x-y < n \le x}\tau_{2k-1}(n)\Big)^k.$$
\end{proposition}
\begin{proof}
    This is a mildly altered Lemma 4.4 in \cite{AGH26} so it can fit more cleanly with the hypercontractive inequality. We have the reordering
    $$\sum_{x-y < n \le x}g(n) = \sum_{n \le x}f(n)\sum_{\sqrt{\frac{x-y}{n}} < m \le \sqrt{\frac{x}{n}}}1.$$
    By Proposition \ref{prop_hyper}, then
    $$\mathbb{E}|\sum_{x-y < n \le x}g(n)|^{2k} \le (\sum_{n \le x}\mu^2(n)\tau_{2k-1}(n)|\sum_{\sqrt{\frac{x-y}{n}} < m \le \sqrt{\frac{x}{n}}}1|^2)^k.$$
    As a result,
    \begin{align*}
        (\lfloor\sqrt{x/n}\rfloor - \lfloor\sqrt{(x-y)/n}\rfloor)^2  
        \ll& \frac{y^2}{nx} + (\lfloor\sqrt{x/n}\rfloor-\lfloor\sqrt{(x-y)/n}\rfloor).
    \end{align*}
    Then we have
    \begin{equation}\label{eqn:expanded_form}
        \mathbb{E}\Big|\sum_{x-y < n \le x}g(n)\Big|^{2k} \ll_k (\sum_{n \le x}\frac{y^2\mu^2(n)\tau_{2k-1}(n)}{nx} + \sum_{n \le x}\tau_{2k-1}(n)\mu^2(n)\big(\lfloor\sqrt{x/n}\rfloor-\lfloor\sqrt{(x-y)/n}\rfloor\big))^k.
    \end{equation}
    The second sum is equal to
    $$\sum_{x-y < m^2n \le x}\mu^2(n)\tau_{2k-1}(n) \le \sum_{x-y < n \le x}\tau_{2k-1}(n).$$
    Substituting this into display (\ref{eqn:expanded_form}) with Proposition \ref{prop:divisor_sums} applied to the weighted sum completes the proof.
\end{proof}
The key number theory result we use is a Parseval identity for Dirichlet series, allowing us to move between partial sums of coefficients and their corresponding Euler product.
\begin{proposition} (Equation (5.26) in Section 5.1 of \cite{montgomery07}) Let $(a_n)_{n=1}^\infty$ be a sequence of complex numbers and $A(s) = \sum_{n=1}^\infty \frac{a_n}{n^s}$ be its associated Dirichlet series, with associated abscissa of convergence $\sigma_c$. Then for any $\sigma > \max\{\sigma_c,0\}$, we have that
\begin{equation*}
    \int_1^\infty \frac{|\sum_{n \le x} a_n|^2}{x^{1+2\sigma}} dx = \frac{1}{2\pi}\int_{-\infty}^\infty \left|\frac{A(\sigma+it)}{\sigma+it}\right|^2 dt. 
\end{equation*} \label{HA1}
\end{proposition} 
This has been used many times in random multiplicative function literature, most strikingly in Harper's work proving Helson's conjecture \cite{aharper201} (among other results). We will use this slightly differently to Caich's application since we need to use different techniques to show that the various events coming from bounding these products almost surely occur. The approach used is similar to that of Gerspach's and Hardy's papers.
\subsection{Probabilistic results}
We also want to exploit the martingale structure in our sums and random Euler products. To do this, we use the same machinery as Caich and Mastrostefano.
\begin{proposition}
    (Doob's $L^r$ and maximal supermartingale inequalities, Theorem 10.9.4 and 10.9.1 in \cite{gut06}) Let $r > 1$. Suppose the real sequence of random variables and $\sigma$-algebras $\{(X_n,\mathcal{F}_n)\}_{0 \le n \le N}$ is a non-negative submartingale which is bounded in $L^r$. Then
    $$\mathbb{E}\big[(\max_{0 \le n \le N}X_n)^r\big] \le \big(\frac{r}{r-1}\big)^r\max_{0 \le n \le N}\left(\mathbb{E}(X_n^r)\right).$$
    Instead, suppose $\{(X_n,\mathcal{F}_n)\}_{0 \le n \le N}$ is a non-negative supermartingale with $\lambda > 0$. In this case,
    $$\lambda \mathbb{P}(\max_{0 \le n \le N}X_n > \lambda) \le \mathbb{E}[X_0].$$
    \label{Lr-Doob}
\end{proposition}
Finally, we give a key inequality first used in Caich's work which is based on a version of the Azuma--Hoeffding inequality. 
\begin{proposition}
    (Caich, Lemma 3.11 in \cite{caich23}) Let $\{(X_n,\mathcal{F}_n\}_{n \le N}$ be a complex sequence of martingale differences. Assume that $X_n$ is bounded almost surely (we suppose there exists real number $b>0$ such that $|X_n| < b$ almost surely). Further, assume that $|X_n|\le S_n$ where $(S_n)_{n \le N}$ is a sequence of real random variables, where for each $n$, $S_n$ is $\mathcal{F}_{n-1}$ measurable. Define the event $\mathcal{E} := \{\sum_{1 \le n \le N}S_n^2 \le T\}$ where $T>0$ is a deterministic constant. Then for any $\varepsilon>0$,
    $$\mathbb{P}\Big(\big\{\big|\sum_{n \le N}X_n\big| > \varepsilon\Big\} \bigcap \{\mathcal{E}\}\Big) \le 2\exp\Big(-\frac{\varepsilon^2}{10T}\Big).$$ \label{Caich_AH_in}
\end{proposition}
\subsection{Random Euler product and multiplicative chaos estimates}
At various points we will require the following expectation result for random Euler products.
\begin{proposition}[Euler Product Result 2, \cite{aharper201}]\label{prop:MomentsEulerProducts}
     Suppose $400 \leq y < z$ for sufficiently large $z$ with $\alpha,h \in \mathbb{R}$ where $\alpha$ is fixed and $\sigma \geqslant -10/\log z$, then 
    \begin{equation*}\begin{split}
        \mathbb{E}  \Big[\prod_{y < p \leq z}\Big|1+\frac{f(p)}{p^{1/2+\sigma+ih}}\Big|^{2\alpha}\Big] 
        = \exp\Big(\sum_{y < p \leq z}\frac{\alpha^2  + (\alpha^2 - \alpha)\cos(2h\log p)}{p^{1+2\sigma}} + O\Big(\frac{1}{\sqrt{y}\log y}\Big)\Big). \end{split}
    \end{equation*} 
\end{proposition}
The most important ingredient of this proof is the following random Euler product estimate from Proposition 5.3 of \cite{AGH26}.
\begin{proposition}\label{prop:upper_uniform_REP}
    Let $x$ be sufficiently large. Suppose that $0 <  \sigma \le \log\log x/\log x$. Uniformly in $x$ and $\sigma$, we have
    \begin{equation}\label{typical_q}\begin{split}
        &\mathbb{E}\Big[\Big(\frac{1}{\log x}\int_{\mathbb{R}} \frac{|F_x(1/2+\sigma+ih)|^2\log^2(|h|+2) }{\sigma^2+h^2} dh\Big)^{1/3}\Big] \ll (\sigma\log^{4/3} x)^{-1/3}+(\log\log x)^{-1/6}.
    \end{split}\end{equation}
\end{proposition}
\section{Reduction to test points and anatomy bounded sets}\label{sec:setup}
\subsection{Reduction to test points}
This proof combines Caich's work \cite{caich23} with the random Euler product analysis coming from \cite{AGH26}. The first step is to rewrite the sums where we split the sum according to the size of the prime factors of $n$. This was used in the work of Lau--Tenenbaum--Wu \cite{LTW13} to improve the efficiency of the inequality in Lemma \ref{prop_hyper}. The first appearance of this argument was in the work of Hal\'{a}sz on random multiplicative functions \cite{halasz83}.
\begin{lemma}
    (Lemma 2.3, \cite{LTW13}) Let $f$ be a Rademacher multiplicative function, and let $A>0$ be some fixed constant. Then there exists $\gamma := \gamma(A) \in (0,1)$ such that 
    \begin{equation}
        x_i := \lfloor e^{i^{\gamma}}\rfloor \quad (i \geqslant 1)
    \end{equation}
    we have that almost surely
    \begin{equation}
        \max_{x_{i-1} \le x \le x_i}\bigg|\sum_{x_{i-1} < n \le x}f(n) \bigg| \ll \frac{\sqrt{x_i}}{\log^Ax_{i}}, \quad (i \geqslant 1).
    \end{equation}
\end{lemma}
From this, we simply need to apply partial summation to show that
\begin{equation}
    \max_{x \in [x_{i-1}, x_i]}\bigg|\sum_{x_{i-1} < n \le x}\frac{f(n)}{\sqrt{n}} \bigg| \ll \frac{\sqrt{x_{i}}}{\sqrt{x_{i-1}}}\frac{1}{\log x_{i-1}} \ll \frac{1}{\log x_{i-1}}. \label{weighted_test_bound}
\end{equation}
For the rest of this paper, we will choose $\gamma$ far smaller than this for reasons which will become more apparent throughout the proof (for our purposes, we assume $\gamma \le 10^{-3}$). For the extended Rademacher case, we need a slightly reworked version of the above.
\begin{lemma}\label{extended_testpoints}
    Let $A>0$ be some fixed constant. Then there exists $\gamma := \gamma(A) \in (0,1)$ such that 
    \begin{equation}
        x_i := \lfloor e^{i^{\gamma}} \rfloor \quad (i \geqslant 1) \label{point_def}
    \end{equation}
    we have that almost surely
    \begin{equation}
        \max_{x_{i-1} \le x \le x_i}\bigg|\sum_{x_{i-1} < n \le x}g(n) \bigg| \ll \frac{\sqrt{x_i}}{\log^A x_{i}}, \quad (i \geqslant 1). \label{test_point}
    \end{equation}
\end{lemma}
\begin{proof}
    This is a minor adjustment to the setup in Basquin \cite{B17} and Lau--Tenenbaum--Wu \cite{LTW13}. We first define a new set of test points between $x_{i-1}$ and $x_i$, named $u_k$ by
    $$u_k = x_{i-1} + \sum_{j=1}^{k}2^{\nu_j}$$
    for some natural numbers $\nu_j$ which are descending in $j$. There are $\le 2\log(x_i)$ such test points between $x_{i-1}$ and $x_i$. Denote $S_g(u) = \sum_{\substack{n \le u}}g(n)$. Then we have
    \begin{align*}
        &\mathbb{P}\Big(\sup_{X_{\ell-1} < x_i \le X_{\ell}}\sup_{x_{i-1} < x \le x_i}\Big|\sum_{\substack{x_{i-1} < n \le x}}g(n)\Big| > 2\sqrt{x_i}\log^{-A}x_i\Big) \\ \ll& \sum_{X_{\ell-1} < x_i \le X_{\ell}}\mathbb{P}\Big(\sup_{x_{i-1} < x \le x_i}\Big|\sum_{\substack{x_{i-1} < n \le x}}g(n)\Big| > 2\sqrt{x_i}\log^{-A}x_i\Big) \\
        \ll& \sum_{X_{\ell-1} < x_i \le X_{\ell}}\mathbb{P}\Big(\sup_{x_{i-1} < u_k \le x_i}\Big|S_g(u_k) - S_g(x_{i-1})\Big| > \sqrt{x_i}\log^{-A}x_i\Big).
    \end{align*}
    We break this difference into chunks. Let $u=x_{i-1}+(l-1)2^m$ and $v=x_{i-1}+l 2^m$ where $l = \sum_{j=1}^k2^{\nu_j-\nu_k} \ge 1$ and $m = \nu_k \ge 0$ such that $l2^m < x_i - x_{i-1}$. Then by the union bound over the $k$ followed by Markov's inequality to the fourth power gives
    \begin{align*}
        \ll& \sum_{X_{\ell-1} < x_i \le X_{\ell}}\frac{(\log x_i)^{4(A+1)}}{x_{i-1}^2}\sum_{\substack{l > 0; m \ge 0 \\ l2^m \le x_i-x_{{i-1}}}}\mathbb{E}[|S_g(v)-S_g(u)|^4] \\
        \ll&\sum_{X_{\ell-1} < x_i \le X_{\ell}}\frac{(\log x_i)^{4(A+1)}}{x_{i-1}^2}\sum_{\substack{l > 0; m \ge 0 \\ l2^m \le x_i-x_{{i-1}}}}\Big(\big(\sum_{u < n \le v}\tau_3(n)\big)^2 + \big(\frac{(v-u)^2\log^3(2v)}{v}\big)^2\Big).
    \end{align*}
    The final line comes from Proposition \ref{prop:short_sum_count} and the elementary inequality $|A+B|^2 \le 2(A^2+B^2)$. Using that
    $$\big(\sum_{u < n \le v}\tau_3(n)\big)^{2} \ll v^{2/3}(\log 2v)^{52/3}(v-u)^{4/3},$$
    then we see
    \begin{align*}
        &\mathbb{P}\Big(\sup_{X_{\ell-1} < x_i \le X_{\ell}}\sup_{x_{i-1} < x \le x_i}\Big|\sum_{\substack{x_{i-1} < n \le x}}g(n)\Big| > 2\sqrt{x_i}\log^{-A}x_i\Big) \\ \ll& \sum_{X_{\ell-1} < x_i \le X_\ell}(\log x_i)^{4(A+1)+52/3}\Big(\frac{x_i-x_{i-1}}{x_{i}}\Big)^{4/3} \ll \sum_{X_{\ell-1} < x_i \le X_\ell}\frac{1}{i^{(4-(68+12A)\gamma)/3}}
    \end{align*}
    which is certainly summable in $\ell$ for $\gamma$ sufficiently small.
\end{proof}
Using this, it is then sufficient to understand both sums at only the $x_i$, which is the objective of the rest of the paper.
\begin{definition}
    \textit{Let $\eta>0$ and $r \geq 2$. Then we almost surely have for the $(x_i)$ defined in equation (\ref{point_def})}
    \begin{equation}
        |\sum_{\substack{n \le x_i}}g(n)| \ll  \sqrt{x_i}(\log\log x_i))^{1/4+\eta+1/2r}, \qquad |\sum_{n \le x_i}\frac{f(n)}{\sqrt{n}}| \ll (\log\log x_i)^{1/4+\eta+1/2r}.\label{key_in}
    \end{equation} \label{key_prop}
\end{definition}
The combination of equations (\ref{test_point}) and (\ref{key_in}) is sufficient to prove the upper bound, since then we will have for any $\eta>0$ and $r \ge 2$
\begin{align*}
    |\sum_{n \le x}g(n)| \le |\sum_{x_{i-1} < n \le x}g(n)| + |\sum_{n \le x_{i-1}}g(n)| \le& \sup_{x_{i-1} < y \le x_i}|\sum_{x_{i-1} < n \le y}g(n)| + |\sum_{n \le x_{i-1}}g(n)| \\
    \le& \frac{\sqrt{x_i}}{\log x_i} + \sqrt{x_i}(\log\log x_i)^{1/4+\eta+1/2r}
\end{align*}
almost surely. Taking $r=\lfloor1/\varepsilon\rfloor + 1 > 2$ and $\eta = \varepsilon/2$ obtains the result. Similarly in the weighted case, we have
\begin{align*}
    |\sum_{n \le x}\frac{f(n)}{\sqrt{n}}| \le |\sum_{x_{i-1} < n \le x}\frac{f(n)}{\sqrt{n}}| + |\sum_{n \le x_{i-1}}\frac{f(n)}{\sqrt{n}}| \le& \sup_{x_{i-1} < y \le x_i}|\sum_{x_{i-1} < n \le y}\frac{f(n)}{\sqrt{n}}| + |\sum_{n \le x_{i-1}}\frac{f(n)}{\sqrt{n}}| \\
    \le& \frac{1}{\log x_i} + (\log\log x_i)^{1/4+\eta+1/2r}
\end{align*}
as required.
\subsection{Slowly varying test points and very smooth contributions}
To prove Key Proposition \ref{key_prop}, we need to use some of definitions used in \cite{dpc26}, \cite{caich23} and \cite{mastrostefano2022almost}. For $\ell \in \mathbb{N}$, we take $X_{\ell} = \exp(2^{\ell^K})$, where $K=\lfloor\frac{25}{\eta}\rfloor$. In particular, we see that each $x_i$ is contained in the interval $(X_{\ell-1},X_{\ell}]$ for some $\ell \in \mathbb{N}$. This makes our task more difficult than what we see in Lau--Tenenbaum--Wu since we will have more test points $x_i$ in the interval $(X_{\ell-1},X_{\ell}]$. This is why we need the more powerful martingale machinery to assist us. We also define the finite sequence of numbers $(y_j)_{j=0}^J$ where we have that
\begin{equation*}
    y_0 = \exp\Big(\frac{2^{(\ell-1)^K}}{\ell^{2K}}\Big), \quad y_j = \exp\big(e^{j/\ell}\log y_0\big).
\end{equation*}
where $J$ is the smallest integer such that $y_J > X_\ell$. Then we have that $J \ll \ell^K$.
We have the following splitting
\begin{align*}
    \sum_{n \le x_i}g(n) =& \sum_{\substack{n \le x_i \\ P(\kappa(n)) \le y_0}}g(n) + \sum_{\substack{n \le x_i \\ P(\kappa(n)) > y_0}}g(n),\\
    \sum_{n \le x_i}\frac{f(n)}{\sqrt{n}} =& \sum_{\substack{n \le x_i \\ P(\kappa(n)) \le y_0}}\frac{f(n)}{\sqrt{n}} + \sum_{\substack{n \le x_i \\ P(\kappa(n)) >y_0}}\frac{f(n)}{\sqrt{n}}. 
\end{align*}
For now, we focus on the sum over the integers with $P(\kappa(n)) \le y_0$ and leave the second part of the sum to be studied in Section \ref{sec:conditional_variance}. Indeed, we have
$$\sum_{\substack{n \le x_i \\ P(\kappa(n)) \le y_0}}g(n) = \sqrt{x}\sum_{P(n) \le y_0}\frac{f(n)}{\sqrt{n}} - \sqrt{x}\sum_{\substack{n > x_i \\ P(n) \le y_0}}\frac{f(n)}{\sqrt{n}} - \sum_{\substack{n \le x_i \\ P(n) \le y_0}}f(n)\Big\{\sqrt{x/n}\Big\}.$$
The additional term involving the fractional part is not present in the weighted Rademacher case, and the sum over the integers with all prime factors less than $y_0$ is an Euler product. Define the events
\begin{align*}
    \mathcal{A}_\ell^{(e)} := \Big\{\sup_{X_{\ell-1} < x_i \le X_\ell}\Big|\frac{1}{\sqrt{x_i}}\sum_{\substack{n \le x_i \\ P(\kappa(n)) \le y_0}}g(n)\Big|>\ell^{K/4}\Big\}, \\
    \mathcal{A}_\ell^{(w)} := \Big\{\sup_{X_{\ell-1} < x_i \le X_\ell}\Big|\sum_{\substack{n \le x_i \\ P(n) \le y_0}}\frac{f(n)}{\sqrt{n}}\Big|>\ell^{K/4}\Big\}.
\end{align*}
\begin{lemma}\label{lemma:smooth_set}
    Both $\mathbb{P}(\mathcal{A}_\ell^{(e)})$ and $\mathbb{P}(\mathcal{A}_\ell^{(w)})$ are summable in $\ell$.
\end{lemma}
\begin{proof}
    By the triangle inequality, we have
    \begin{align*}
        \mathbb{P}(\mathcal{A}_\ell^{(e)}) \le \mathbb{P}(|F_{y_0}(1/2)| > \ell^{K/4}/3) &+ \mathbb{P}(\sup_{X_{\ell-1} < x_i \le X_{\ell}}|\sum_{\substack{n > x_i \\ P(n) \le y_0}}\frac{f(n)}{\sqrt{n}}| > \ell^{K/4}/3)\\ + \mathbb{P}(\sup_{X_{\ell-1} < x_i \le X_{\ell}}&|\frac{1}{\sqrt{x_i}}\sum_{\substack{n \le x_i \\ P(n) \le y_0}}f(n)\{\sqrt{x/n}\}| > \ell^{K/4}/3).
    \end{align*}
    The corresponding weighted version only involves the first two sums. Then by union bound and Chebychev's inequality, the final probability is bounded by
    $$\mathbb{E}|\sum_{\substack{n \le x_i \\ P(n) \le y_0}}f(n)\{\sqrt{x/n}\}|^2 \ll \Psi(x_i,y_0) \ll x_i \log y_0 \exp(-u\log u(1+o(1)))$$
    where $u=\frac{\log x_i}{\log y_0}$ using Theorem 7.6 in \cite{montgomery07}. Using that $u \ge \ell^{2K}$ then
    $$\frac{1}{\ell^{K/2}x_i}\sum_{X_{\ell-1} < x_i \le X_\ell}\mathbb{E}|\sum_{\substack{n \le x_i \\ P(n) \le y_0}}f(n)\{\sqrt{x/n}\}|^2 \ll \exp(\frac{\ell^K\log 2}{\gamma}-C\ell^{2K})$$
    where $C>0$, which is summable in $\ell$ given $\ell$ is sufficiently large. The middle probability can once again be handled by the union bound and Chebychev's inequality, and we have this is
    \begin{equation}\label{eq:preRankin}
        \ll \frac{1}{\ell^{K/2}}\sum_{X_{\ell-1} < x_i \le X_\ell}\sum_{\substack{n > x_i \\ P(n) \le y_0}}\frac{1}{n}.
    \end{equation}
    By Rankin's trick, we get
    $$\sum_{\substack{n > x_i \\ P(n) \le y_{0}}}\frac{1}{n} \le x_i^{-1/\log y_{0}}\prod_{p \le y_{0}}\bigg(1-\frac{1}{p^{1-\frac{1}{\log y_{0}}}}\bigg)^{-1} \ll \log y_{0}\exp(-C\ell^{2K}).$$
    Placing this estimate back into equation (\ref{eq:preRankin}), then this is
    $$\ll \ell^{-K/2}\exp(\log 2((1+\frac{1}{\gamma})\ell^K -\ell^{2K}))$$
    which is again summable in $\ell$. As for the final term, we can simply take Markov to the $q=1/2$, and we have
    $$\mathbb{P}(|F_{y_0}(1/2)| > \ell^{K/4}/3) \ll \ell^{-K/8}\mathbb{E}|F_{y_0}(1/2)|^{1/2} \ll \ell^{-K/8}(\log y_0)^{-1/8}$$
    using Proposition \ref{prop:MomentsEulerProducts} as this quantity is independent of the index $i$, so there is no union bound required. This is also summable in $\ell$ (for sufficiently small $\varepsilon$), which concludes the proof.
\end{proof}
\section{Conditional variance}\label{sec:conditional_variance}
\subsection{Introducing an integral smoothing}
We make the final reductions to this sum. We have for both weighted and extended Rademacher cases
$$\sum_{\substack{n \le x_i \\ P(\kappa(n)) > y_0}}g(n) = \sum_{y_0 < p \le x_i}f(p)\sum_{\substack{n \le x_i/p \\ P(\kappa(n)) < p}}g(n), \quad \sum_{\substack{n \le x_i \\ P(n) > y_0}}\frac{f(n)}{\sqrt{n}} =  \sum_{y_0 < p \le x_i}\frac{f(p)}{\sqrt{p}}\sum_{\substack{n \le x_i/p \\ P(n) < p}}\frac{f(n)}{\sqrt{n}}.$$
Using the identification
$$A_p^{(e)}(x) = f(p)\sum_{\substack{n \le x/p \\ P(\kappa(n)) < p}}g(n), \qquad A_p^{(w)}(x) = \frac{f(p)}{\sqrt{p}}\sum_{\substack{n \le x/p \\ P(n) < p}}\frac{f(n)}{\sqrt{n}},$$
the sums $\sum_p A_p^{(e)}(x)$ and $\sum_{p}A_p^{(w)}(x)$ are martingale difference sequences in $p$. Our intention is to apply Proposition \ref{Caich_AH_in}. After the previous estimates, it remains to control in the extended case
$$V_{\ell}^{(e)}(x_i) = \sum_{y_0 < p \le x_i}\mathbb{E}\Big[\big|A_p^{(e)}(x_i)\big|^2|\mathcal{F}_{p^-}\Big]=\sum_{y_{0} < p \le x_i}|\sum_{\substack{ n\le x_i/p \\ P(\kappa(n)) < p}}g(n)|^2.$$
We will follow the treatment by Caich to reduce this to a random Euler product involving extended Rademacher $f(p)$, which can then be transferred to a corresponding statement about an integral of a standard Rademacher random Euler product. This reduction is identical in the weighted squarefree Rademacher case, so outside of summarising the final quantities we omit the deduction. \par
Indeed, let $\mathcal{X}$ be a large quantity defined in terms of $x_i$ to be determined later. Then we apply the triangle inequality to see
\begin{align*}
    V_{\ell}^{(e)}(x_i) \ll& \sum_{\substack{j\ge 1 \\ y_{j-1} \le x_i}} \sum_{y_{j-1} < p \le y_j}\frac{\mathcal{X}}{p}\int_{p}^{p(1+1/\mathcal{X})}\Big|\sum_{m \ge 1}\sum_{\substack{n \le x/tm^2 \\ P(n) < p}}f(n)\Big|^2dt \\
    &+ \sum_{y_{0} < p \le x_i}\frac{\mathcal{X}}{p}\int_{p}^{p(1+1/\mathcal{X})}\Big|\sum_{m \ge 1}\sum_{\substack{x/tm^2 < n \le x/pm^2 \\ P(n) < p}}f(n)\Big|^2dt \\
    &=: \sum_{\substack{j\ge 1 \\ y_{j-1} \le x_i}}\mathcal{C}_{j,\ell}^{(e)}(x_i,y_j) + \mathcal{D}^{(e)}_\ell(x_i).
\end{align*}
We will continue editing the first integral. Then we have that we exchange the order of summation and integration, and find
$$\mathcal{C}_{j,\ell}^{(e)}(x_i) \ll x_i\int_{x/y_{j}(1+1/\mathcal{X})}^{x/y_{j-1}}\Big|\sum_{m \ge 1}\sum_{\substack{n \le z/m^2 \\ P(n) < p}}f(n)\Big|^2\sum_{ p \in (y_{j-1},y_j] \cap (\frac{x_i}{z(1+1/\mathcal{X})}, \frac{x_i}{z}]}\frac{\mathcal{X}}{p}\frac{dz}{z^2}$$
using the substitution $z = x_i/t$. The sum over $p$ is handled by the prime number theorem, and we see that
$$\sum_{\frac{x_i}{z(1+1/\mathcal{X})} < p \le \frac{x_i}{z}}\frac{\mathcal{X}}{p} \ll \frac{1}{\log y_{j} }, \qquad \sum_{\frac{x_i}{z(1+1/\mathcal{X})} < p \le \frac{x_i}{z}}\frac{\mathcal{X}}{p^2} \ll \frac{z}{x_i\log y_j}.$$
Consequently, we can bound $\mathcal{C}_{j,\ell}^{(e)}(x_i)$ by 
$$\frac{x_i}{\log y_j}\int_{x_i/y_{j}(1+1/\mathcal{X})}^{x_i/y_{j-1}}\sup_{y_{j-1} < p \le y_j}\Big|\sum_{m \ge 1}\sum_{\substack{n \le z/m^2 \\ P(n) < p }}f(n)\Big|^2\frac{dz}{z^{2+2\sigma_\ell}}$$
where $\sigma_\ell = 1/\log X_\ell$ which only shifts things by a fixed multiplicative constant since $x_i \le X_\ell$ by definition. Finally, we set
\begin{align*}
    Q_{j,\ell}^{(e)}(x_i) :=& \frac{1}{\log y_j }\int_{x_i/y_j(1+1/\mathcal{X})}^{x_i/y_{j-1}}\sup_{y_{j-1}<p \le y_j}\Big|\sum_{m \ge 1}\sum_{\substack{n \le z/m^2 \\ P(n) < p}}f(n)\Big|^2\frac{dz}{z^{2+2\sigma_\ell}}.
\end{align*}
This leads us to defining the following exceptional sets with $T(\ell) = \ell^{10}$
\begin{align*}
    \mathcal{Q}_{\ell,1}^{(e)} :=& \Big\{\sup_{X_{\ell-1} < x_i \le X_\ell}\sup_{\substack{1 \le j \le J \\ 1 \le \frac{\log x_i}{\log y_{j-1}} \le \ell^{2K}}}Q_{j,\ell}^{(e)}(x_i) > \frac{T(\ell)}{\ell^{K(1-2/r)/2}\ell\log \ell}\Big\}, \\
    \mathcal{Q}_{\ell,2}^{(e)} :=& \Big\{\sup_{X_{\ell-1} < x_i \le X_\ell}\sup_{\substack{1 \le j \le J \\ \frac{\log x_i}{\log y_{j-1}} > \ell^{2K}}}Q_{j,\ell}^{(e)}(x_i) > \ell^{-3K/2}\Big\}, \\
    \mathcal{D}_\ell^{(e)} :=& \Big\{\sup_{X_{\ell-1} < x_i \le X_\ell}D_\ell^{(e)}(x_i) > \frac{x_i}{\ell^{K/2}}\Big\}.
\end{align*}
We have in the weighted case
\begin{equation*}\begin{split}
    Q_{j,\ell}^{(w)}(x_i) :=& \frac{1}{\log y_j }\int_{x_i/y_j(1+1/\mathcal{X})}^{x_i/y_{j-1}}\sup_{y_{j-1}<p \le y_j}\Big|\sum_{\substack{n \le z \\ P(n) < p}}\frac{f(n)}{\sqrt{n}}\Big|^2\frac{dz}{z^{1+2\sigma_\ell}}, \\ 
    \mathcal{D}_\ell^{(w)}(x_i) :=&\sum_{y_{0} < p \le x_i}\frac{\mathcal{X}}{p^2}\int_{p}^{p(1+1/\mathcal{X})}\Big|\sum_{\substack{x/t < n \le x/p \\ P(n) < p}}\frac{f(n)}{\sqrt{n}}\Big|^2dt,
\end{split}
\end{equation*}
and $\mathcal{Q}_{\ell,k}^{(w)}$ are defined in the same way as their extended counterparts, only with their $(w)$ versions used in their definitions. Instead, we have for the weighted case
$$\mathcal{D}_\ell^{(w)} := \Big\{\sup_{X_{\ell-1} < x_i \le X_\ell}D_\ell^{(w)}(x_i) > \ell^{-K/2}\Big\}.$$
We introduce this splitting in the $y_j$ since we expect only the largest $y_j$ to contribute to the large fluctuations, and there are $\ll \ell\log \ell$ of such $y_j$. This and the integral Minkowski argument from \cite{dpc26} allow us to achieve the $1/4$ exponent in our theorems.
\subsection{Short interval and very smooth exceptional sets}
To handle the event $\mathcal{D}_\ell^{(e)}$ and $\mathcal{D}_\ell^{(w)}$, we will take very high moments of the quantity and use our choice of $\mathcal{X}$ to dominate things accordingly.
\begin{lemma}\label{lem:short_extended_excep}
    We have that $\mathbb{P}(\mathcal{D}_\ell^{(e)})$ is summable in $\ell$.
\end{lemma}
\begin{proof}
Suppose $q \ge 1$ is a fixed parameter chosen later. We consider
$$\bigg(\mathbb{E}\bigg[\Big(\sum_{y_0 < p \le x_i }\frac{\mathcal{X}}{p}\int_{p}^{p(1+1/\mathcal{X})}\Big|\sum_{m \ge 1}\sum_{\substack{x_i/tm^2 < n \le x_i/pm^2 \\ P(n) < p}}f(n)\Big|^2dt\Big)^q\bigg]\bigg)^{1/q}.$$
By Minkowski inequality, this is
$$\sum_{y_0 < p \le x_i}\bigg(\mathbb{E}\bigg[\Big(\frac{\mathcal{X}}{p}\int_{p}^{p(1+1/\mathcal{X})}\Big|\sum_{m \ge 1}\sum_{\substack{x_i/tm^2 < n \le x_i/pm^2 \\ P(n) < p}}f(n)\Big|^2dt\Big)^q\bigg]\bigg)^{1/q}.$$
The integral is appropriately normalised, so we can apply H\"older's inequality to the $q$-th power and obtain the further upper bound
\begin{equation}\label{eq:hyper_manipulation}
    \sum_{y_0 < p \le x_i}\bigg(\frac{\mathcal{X}}{p}\int_{p}^{p(1+1/\mathcal{X})}\mathbb{E}\bigg[\Big|\sum_{m \ge 1}\sum_{\substack{x_i/tm^2 < n \le x_i/pm^2 \\ P(n) < p}}f(n)\Big|^{2q}\bigg]dt\bigg)^{1/q}.
\end{equation}
We have to proceed slightly differently because of the square contributions. Here, we immediately apply Proposition \ref{prop:short_sum_count}, and find this is
$$\ll \sum_{y_0 < p \le x_i}\Big(\frac{\mathcal{X}}{p}\int_{p}^{p(1+1/\mathcal{X})}\Big(\frac{(x_i/p-x_i/t)^2\log^{2q-1}(2x_i/p)}{x_i/p}+\sum_{x_i/t < n \le x_i/p}\tau_{2q-1}(n)\Big)^qdt\Big)^{1/q}.$$
We have removed the smoothness restriction here. Then $|A+B|^{2q} \le 2^{2q}(|A|^{2q}+|B|^{2q})$, so
$$\ll_q \sum_{y_0 < p \le x_i}\Big(\frac{\mathcal{X}}{p}\int_{p}^{p(1+1/\mathcal{X})}\Big(\frac{(x_i/p-x_i/t)^2\log^{2q-1}(2x_i/p)}{x_i/p}\Big)^q +\Big(\sum_{x_i/t < n \le x_i/p}\tau_{2q-1}(n)\Big)^qdt\Big)^{1/q}.$$
To split these terms, we use subadditivity as $1/q \le 1$. For the first term, we use a length max argument, and see this is
$$\Big(\frac{\mathcal{X}}{p}\int_{p}^{p(1+1/\mathcal{X})}\Big(\frac{(x_i/p-x_i/t)^2\log^{2q-1}(2x_i/p)}{x_i/p}\Big)^qdt\Big)^{1/q} \ll_q \frac{x_i\log^{2q-1}(2x_i)}{p\mathcal{X}^2}.$$
The second term, we instead undo the integral smoothing to find it is
$$\ll \sum_{y_0 < p \le x_i}\sum_{\frac{x_i}{p(1+1/\mathcal{X})} < n \le \frac{x_i}{p}}\tau_{2q-1}(n).$$
Then using the convolution identity $1*\tau_{2q-1} = \tau_{2q}$, we have
$$\ll \sum_{y_0 < p \le x_i}\sum_{\frac{x_i}{p(1+1/\mathcal{X})} < n \le \frac{x_i}{p}}\tau_{2q-1}(n) \le \sum_{\frac{x_i}{1+1/\mathcal{X}} < n \le x_i}\sum_{p|n}\tau_{2q-1}(n/p) \le \sum_{\frac{x_i}{1+1/\mathcal{X}} < n \le x_i}\tau_{2q}(n).$$
Finally we can use the Shiu bound from Proposition \ref{prop:divisor_sums} with $h \ll x_i/\mathcal{X}$ (provided $\mathcal{X} = o(\sqrt{x_i})$), and we get the above is $\ll \tfrac{x_i\log^{2q}x_i}{\mathcal{X}}$
Putting these together, we have
\begin{align*}
    &\bigg(\mathbb{E}\bigg[\Big(\sum_{y_0 < p \le x_i }\frac{\mathcal{X}}{p}\int_{p}^{p(1+1/\mathcal{X})}\Big|\sum_{m \ge 1}\sum_{\substack{x_i/tm^2 < n \le x_i/pm^2 \\ P(n) < p}}f(n)\Big|^2dt\Big)^q\bigg]\bigg)^{1/q} \\
    \ll_q& \frac{x_i\log^{2q}(2x_i)}{\mathcal{X}} + \sum_{y_0 < p \le x_i}\frac{x_i\log^{2q-1}(2x_i)}{p\mathcal{X}^2} \ll \frac{x_i\log^{2q}(2x_i)}{\mathcal{X}}.
\end{align*}
By the union bound and Markov's inequality to the $q$ power, we see
\begin{align*}
    \mathbb{P}(\mathcal{D}_\ell) \ll \ell^{Kq/2}x_i^{-q}\sum_{X_{\ell-1} < x_i \le X_\ell}\bigg(\frac{x_i\log^{2q} (2x_i)}{\mathcal{X}} \bigg)^q.
\end{align*}
Choosing $\mathcal{X} = \log^{2q+1}(2x_i)$ shows the above is
$$\ll \sum_{X_{\ell-1} < x_i \le X_\ell}\Big(\frac{\ell^K}{\log x_i}\Big)^q$$
which is summable for the choice of integer $q > 1/\gamma$.
\end{proof}
The weighted version of this lemma is slightly different to the above and involves some non trivial adjustments.
\begin{lemma}\label{lem:short_weighted_excep}
    We have that $\mathbb{P}(\mathcal{D}_\ell^{(w)})$ is summable in $\ell$.
\end{lemma}
\begin{proof}
Let $r \geqslant 1$ be a constant to be chosen later. Then we can apply the union bound and Markov's inequality to the $r$-power followed by Minkowski's inequality to see
\begin{align*}
    \mathbb{P}(\mathcal{D}_\ell^{(w)})
    \ll& \sum_{X_{\ell-1} < x_i \le X_\ell}\ell^{3Kr/2}((\mathbb{E}[({D}_\ell^{(w)}(x_i))^r])^{1/r})^r.
\end{align*}
We next proceed with finding suitable upper bounds for $(\mathbb{E}[(D_\ell^{(w)}(x_i))^r])^{1/r}$. We follow the approach of Harper \cite{harper202} and Hardy \cite{shardy23}. We have
$$(\mathbb{E}[(D_\ell^{(w)}(x_i))^r])^{1/r} \le \sum_{y_{0} < p \le x_i}\frac{1}{p}\big(\mathbb{E}\Big(\frac{\mathcal{X}}{p}\int_{p}^{p(1+1/\mathcal{X})}\big|\sum_{\substack{x_i/t < n \le x_i/p\\P(n) <p}}\frac{f(n)}{\sqrt{n}}\big|^2 dt\Big)^r\big)^{1/r}.$$
Since we have normalised the integral above, we can apply H\"{o}lder's inequality,
\begin{align*}
    (\mathbb{E}[(D_\ell^{(w)}(x_i))^r])^{1/r} \le \sum_{y_{0} < p \le x_i}\frac{1}{p}\Big(\frac{\mathcal{X}}{p}\int_{p}^{p(1+1/\mathcal{X})}\mathbb{E}\big|\sum_{\substack{x_i/t < n \le x_i/p \\ P(n) <p}}\frac{f(n)}{\sqrt{n}}\Big|^{2r} dt\big)^{1/r}.
\end{align*}
We split the above integral according to $p > \frac{x_i}{1+\mathcal{X}}$. In the first case, we find that $\frac{x_i}{p}-\frac{x_i}{p(1+1/\mathcal{X})} < 1$, so the inner sum has at most one term. This gives that for this range
$$\frac{\mathcal{X}}{p}\int_{p}^{p(1+1/\mathcal{X})}\mathbb{E}\Big|\sum_{\substack{x_i/t < n \le x_i/p\\P(n) < p}}\frac{f(n)}{\sqrt{n}}\Big|^{2r}dt \ll \frac{\mathcal{X}}{p}\int_{p}^{p(1+1/\mathcal{X})}\Big(\frac{t}{x_i}\Big)^{r}dt \ll \Big(\frac{p}{x_i}\Big)^r.$$
We take the $1/r$ power of this before we can sum this over the primes, so the $p$ powers cancel, with the $x_i$ remaining. Summing over the $p>\frac{x_i}{1+\mathcal{X}}$ and using Abel summation with the prime number theorem we have
$$\frac{1}{x_i}\sum_{\frac{x_i}{1+\mathcal{X}} < p \le x_i}1 \ll \frac{1}{\log x_i}.$$
For the range $p \le \frac{x_i}{1+\mathcal{X}}$, then we need to proceed slightly differently. We follow Hardy's treatment of this problem \cite{shardy23}, where we use Proposition \ref{prop_hyper} followed by Cauchy--Schwarz. For the first step, we have
$$\mathbb{E}\Big|\sum_{\substack{x_i/t < n \le x_i/p\\P(n) < p}}\frac{f(n)}{\sqrt{n}}\Big|^{2r} \le \Big(\sum_{\substack{x_i/t < n \le x_i/p\\P(n) <p}}\frac{\tau_{2r-1}(n)}{n}\Big)^{r}.$$
Applying Cauchy--Schwarz, we find
$$\Big(\sum_{\substack{x_i/t < n \le x_i/p\\P(n) <p}}\frac{\tau_{2r-1}(n)}{n}\Big)^{r} \le \Big(\sum_{\substack{x_i/t < n \le x_i/p\\P(n)<p}}\frac{1}{n^2}\Big)^{r/2}\Big(\sum_{\substack{x_i/t < n \le x_i/p\\P(n) < p}}\left(\tau_{2r-1}(n)\right)^2\Big)^{r/2}.$$
We have
$$\Big(\sum_{\substack{x_i/t < n \le x_i/p\\P(n) <p}}\frac{1}{n^2}\Big)^{r/2} \ll \Big(\sum_{\frac{x_i}{p(1+1/X)} < n \le x_i/p}\frac{1}{n^2}\Big)^{r/2} \ll \Big(\frac{p}{\mathcal{X} x_i}\Big)^{r/2}.$$
The second term we apply Cauchy-Schwarz to this sum and observe that $(\tau_{2r-1}(n))^{2} \le \tau_{4r^2-4r+1}(n)$ (as seen in the paper of Benatar--Nishry--Rodgers \cite{benatar22} equations 3.2 and 3.3),
$$\Big(\sum_{\substack{x_i/t < n \le x_i/p\\P(n) < p}}(\tau_{2r-1}(n))^2\Big)^{r/2} \le \Big(\sum_{\substack{x_i/t < n \le x_i/p\\P(n) < p}}\tau_{4r^2-4r+1}(n)\Big)^{r/2} \ll \Big(\frac{x_i\log^{4r^2-4r}x_i}{p}\Big)^{r/2}.$$
Putting this all together, we have that for $p \le \frac{x_i}{1+\mathcal{X}}$
$$\mathbb{E}\Big|\sum_{\substack{x_i/t < n \le x_i/p\\P(n)<p}}\frac{f(n)}{\sqrt{n}}\Big|^{2r} \ll \Big(\frac{\log^{4r^2-4r}x_i}{\mathcal{X}}\Big)^{r/2}.$$
This leads to us having
$$(\mathbb{E}[(D_\ell^{(w)}(x_i))^r])^{1/r} \ll \sum_{p \le \frac{x_i}{1+\mathcal{X}}}\frac{1}{p}\frac{ \log^{2r^2-2r} x_i} {\sqrt{\mathcal{X}}} + \frac{1}{\log x_i}.$$
This can be substituted into the initial bound on $\mathbb{P}(\mathcal{D}_\ell^{(w)})$;
$$\mathbb{P}(\mathcal{D}_\ell^{(w)}) \ll \sum_{X_{\ell-1} < x_i \le X_\ell}\Big(\sum_{y_0 < p \le \frac{x_i}{1+\mathcal{X}}}\frac{1}{p}\frac{\ell^{3K/2}\log^{2r^2-2r}x_i}{\sqrt{\mathcal{X}}}\Big)^{r} + \sum_{X_{\ell-1} < x_i \le X_\ell}\Big(\frac{\ell^{3K/2}}{\log x_i}\Big)^r.$$
We then choose $\mathcal{X} = \log^{4r^2-4r+2}x_i$. With this choice, we have
$$\mathbb{P}(\mathcal{D}_\ell^{(w)}) \ll \sum_{X_{\ell-1} < x_i \le X_\ell} \Big(\frac{\ell^{3K/2}}{\log x_i}\Big)^r.$$
Finally, we choose $r > 1/\gamma$ where $\gamma$ is provided in Section \ref{sec:setup}, which is sufficient.
\end{proof}
We also handle the summability of $\mathbb{P}(\mathcal{Q}_{\ell,2}^{(e)})$ and its weighted counterpart.
\begin{lemma}\label{lem:smooth_extended_sum}
    Both $\mathbb{P}(\mathcal{Q}_{\ell,2}^{(e)})$ and $\mathbb{P}(\mathcal{Q}_{\ell,2}^{(w)})$ are summable in $\ell$.
\end{lemma}
\begin{proof}
    We have to augment this event so we can manipulate it more effectively. Recall we are considering the probability
    $$\mathbb{P}\Big(\sup_{X_{\ell-1} < x_i \le X_{\ell}}\sup_{\substack{1 \le j \le J \\ \frac{\log x_i}{\log y_{j-1}} > \ell^{2K}}}\frac{1}{\log y_{j}} \int_{x_i/y_j(1+1/\mathcal{X})}^{x_i/y_{j-1}} \sup_{y_{j-1} < p \le y_j}\Big|\sum_{m \ge 1}\sum_{\substack{n \le z/m^2 \\ P(n) < p}}f(n)\Big|^2\frac{dz}{z^{2+2\sigma_\ell}} > \ell^{-3K/2}\Big).$$
    Using the same splitting device we established in Lemma \ref{lemma:smooth_set}, then this is
    \begin{align*}
        \le& \mathbb{P}\Big(\sup_{X_{\ell-1} < x_i \le X_{\ell}}\sup_{\substack{1 \le j \le J \\ \frac{\log x_i}{\log y_{j-1}} > \ell^{2K}}}\frac{1}{\log y_{j}} \int_{x_i/y_j(1+1/\mathcal{X})}^{x_i/y_{j-1}} \sup_{y_{j-1} < p \le y_j}\Big|\sum_{\substack{P(n) < p}}\frac{f(n)}{\sqrt{n}}\Big|^2\frac{dz}{z^{1+2\sigma_\ell}} > \frac{1}{9\ell^{3K/2}}\Big) \\
         +& \mathbb{P}\Big(\sup_{X_{\ell-1} < x_i \le X_{\ell}}\sup_{\substack{1 \le j \le J \\ \frac{\log x_i}{\log y_{j-1}} > \ell^{2K}}}\frac{1}{\log y_{j}} \int_{x_i/y_j(1+1/\mathcal{X})}^{x_i/y_{j-1}} \sup_{y_{j-1} < p \le y_j}\Big|\sum_{\substack{n > z \\ P(n) < p}}\frac{f(n)}{\sqrt{n}}\Big|^2\frac{dz}{z^{1+2\sigma_\ell}} > \frac{1}{9\ell^{3K/2}}\Big) \\
         +& \mathbb{P}\Big(\sup_{X_{\ell-1} < x_i \le X_{\ell}}\sup_{\substack{1 \le j \le J \\ \frac{\log x_i}{\log y_{j-1}} > \ell^{2K}}}\frac{1}{\log y_{j}} \int_{x_i/y_j(1+1/\mathcal{X})}^{x_i/y_{j-1}} \sup_{y_{j-1} < p \le y_j}\Big|\sum_{\substack{n \le z \\ P(n) < p}}f(n)\{\sqrt{z/n}\}\Big|^2\frac{dz}{z^{2+2\sigma_\ell}} > \frac{1}{9\ell^{3K/2}}\Big). 
    \end{align*}
    The first two quantities appear in both the weighted and extended case, with the last probability only being considered in the extended case. We can evaluate the integral in the first probability now, and the integral is $\ll \log(y_j/y_{j-1}) \ll \log y_{j}$, so this probability is now independent of the $x_i$ once we extend the supremum condition over all $j \le J$. Then for some absolute constant $C>6$, we are left to determine
    $$\mathbb{P}(\sup_{1 \le j \le J} \sup_{y_{j-1} < p \le y_j}|F_p(1/2)|^2 > \ell^{-3K/2}/C).$$
    In order to collapse the suprema, we see $F_p(1/2)$ is real valued, and is a product of positive terms, so it is positive and $F_p(1/2) = |F_p(1/2)|$. We also have for $Y_p = \sqrt{F_p(1/2)}$
    $$\mathbb{E}\Big[Y_p|\mathcal{F}_{p^-}\Big] = \mathbb{E}\Big[(1+\frac{f(p)}{p^{1/2}})^{1/2}\Big]Y_{p^-} \le (1-\frac{1}{8p}+O(p^{-2}))Y_{p^-} \le Y_{p^-},$$
    so is a non negative supermartingale for sufficiently large $p$. We then have
    $$\mathbb{P}(\sup_{1 \le j \le J}\sup_{y_{j-1} < p \le y_j}|F_p(1/2)|^2 > \ell^{-3K/2}/C) \le \mathbb{P}(\sup_{y_0 \le p < y_J}Y_p > \ell^{{-3K/8}}/C')$$
    for $C' > C$ an absolute constant say. Then by the maximal supermartingale inequality in Proposition \ref{Lr-Doob}, we have
    $$\mathbb{P}(\sup_{y_0 \le p < y_J}Y_p > \ell^{{-3K/8}}/C') \ll \ell^{3K/8}\mathbb{E}[Y_{y_0}] \ll \frac{\ell^{3K/8}}{(\log y_0)^{1/8}}$$
    which is summable in $\ell$.\\
    For the second probability, we apply Markov's inequality to the first power and the union bound, so
    \begin{align*}
        &\mathbb{P}(\sup_{X_{\ell-1} < x_i \le X_{\ell}}\sup_{\substack{1 \le j \le J \\ \frac{\log x_i}{\log y_{j-1}} > \ell^{2K}}}\frac{1}{\log y_{j}} \int_{x_i/y_j(1+1/\mathcal{X})}^{x_i/y_{j-1}} \sup_{y_{j-1} < p \le y_j}\Big|\sum_{\substack{n > z \\ P(n) < p}}\frac{f(n)}{\sqrt{n}}\Big|^2\frac{dz}{z^{1+2\sigma_\ell}} > \frac{1}{9\ell^{3K/2}}) \\
        \ll&  \ell^{3K/2}\sum_{X_{\ell-1} < x_i \le X_{\ell}}\sum_{\substack{1 \le j \le J \\ \frac{\log x_i}{\log y_{j-1}} > \ell^{2K}}}\frac{1}{\log y_{j}} \int_{x_i/y_j(1+1/\mathcal{X})}^{x_i/y_{j-1}} \mathbb{E}\big[\sup_{y_{j-1} < p \le y_j}\Big|\sum_{\substack{n > z \\ P(n) < p}}\frac{f(n)}{\sqrt{n}}\Big|^2\big]\frac{dz}{z^{1+2\sigma_\ell}}.
    \end{align*}
    We have $\sum_{\substack{n > z \\ P(n) < p}}\frac{f(n)}{\sqrt{n}}$ is a martingale with respect to the filtration $\mathcal{F}_{p^-}$, so the mean square of the sum is a submartingale. This allows us to apply Proposition \ref{Lr-Doob} with $r=2$, so followed by Proposition \ref{prop_hyper}, this reduces to bounding
    $$\ell^{3K/2}\sum_{X_{\ell-1} < x_i \le X_{\ell}}\sum_{\substack{1 \le j \le J \\ \frac{\log x_i}{\log y_{j-1}} > \ell^{2K}}}\frac{1}{\log y_{j}} \int_{x_i/y_j(1+1/\mathcal{X})}^{x_i/y_{j-1}} \sum_{\substack{n > z \\ P(n) < y_j}}\frac{1}{n}\frac{dz}{z^{1+2\sigma_\ell}}.$$
    The sum is dealt with using Rankin's trick, and evaluates to
    $$\sum_{\substack{n > z \\ P(n) < y_j}}\frac{1}{n} \ll \log y_j \exp(-\frac{\log z}{\log y_j}) \ll \log y_j \exp(-\ell^{2K}/2).$$
    Substituting the bound back into the integral, we find that this quantity is
    $$\ll \ell^{5K/2}\exp(\ell^K/\gamma - \ell^{2K}/2)$$
    which is certainly summable in $\ell$. \par
    As for the third probability, we can immediately apply Markov's inequality to the first power, and use that $\sum_{\substack{n \le z \\ P(n) < p}}f(n)\{\sqrt{z/n}\}$ is also a martingale with respect to the filtration $\mathcal{F}_{p^-}$, and consequently its mean square is a submartingale with respect to the same filtration. Applying Proposition \ref{Lr-Doob} at $r=2$ gives this is
    \begin{align*}
        &\mathbb{P}(\sup_{X_{\ell-1} < x_i \le X_{\ell}}\sup_{\substack{1\le j \le J \\ \frac{\log x_i}{\log y_{j-1}} > \ell^{2K}}}\frac{1}{\log y_{j}} \int_{x_i/y_j(1+1/\mathcal{X})}^{x_i/y_{j-1}} \sup_{y_{j-1} < p \le y_j}\Big|\sum_{\substack{n \le z \\ P(n) < p}}f(n)\{\sqrt{z/n}\}\Big|^2\frac{dz}{z^{2+2\sigma_\ell}} > \frac{1}{9\ell^{3K/2}}) \\
        \ll& \ell^{3K/2}\sum_{X_{\ell-1} < x_i \le X_{\ell}}\sum_{\substack{1 \le j \le J \\ \frac{\log x_i}{\log y_{j-1}} > \ell^{2K}}}\frac{1}{\log y_{j}} \int_{x_i/y_j(1+1/\mathcal{X})}^{x_i/y_{j-1}} \mathbb{E}\Big|\sum_{\substack{n \le z \\ P(n) < y_j}}f(n)\{\sqrt{z/n}\}\Big|^2\frac{dz}{z^{2+2\sigma_\ell}}.
    \end{align*}
    Proposition \ref{prop_hyper} followed by the coefficient bound $|\{\sqrt{z/n}\}|^2 \le 1$, and the smooth number estimate $\Psi(z,y_j) \ll z\log y_j \exp(-\ell^{2K}/2)$ shows the above is
    \begin{align*}
        \ll& \ell^{3K/2}\sum_{X_{\ell-1} < x_i \le X_{\ell}}\sum_{\substack{1 \le j \le J \\ \frac{\log x_i}{\log y_{j-1}} > \ell^{2K}}}\frac{\exp(-\ell^{2K}/2)\log y_j}{\log y_{j}} \int_{x_i/y_j(1+1/\mathcal{X})}^{x_i/y_{j-1}} \frac{dz}{z^{1+2\sigma_\ell}} \\
        \ll& \ell^{5K/2}\exp(\ell^K/\gamma-\ell^{2K}/2)
    \end{align*}
    which is again summable in $\ell$.
\end{proof}
Once again, the weighted version of this proof is even easier, Rankin's trick is even easier to access as there is no need to handle the terms with square factors.
\subsection{Multiplicative chaos and restricted smooth sums}
To deal with the main contribution to the exceptional sets, we look to Proposition \ref{prop:upper_uniform_REP}. However, to apply the estimate effectively, we set up the following arguments which first appeared in \cite{dpc26}. The main point is that in order to lose the dependence on $i$ so we can use Doob's maximal inequality in an effective manner (as we want a supermartingale in only $j$). Let
\begin{equation}\label{eqn:U_j_def}
    U_j:=N(\sigma_\ell,y_{j})\int_1^\infty\sup_{y_{j-1} \le p \le y_j}\Big|\sum_{m \ge 1}\sum_{\substack{n \le z/m^2 \\ P(n) < p}}f(n)\Big|^2\frac{dz}{z^{2+2\sigma_\ell}}, \qquad N(\sigma,y) = \prod_{p \le y}(1+p^{-1-2\sigma})^{-1}.
\end{equation}
Note $N(\sigma,y) \asymp \tfrac{1}{\log y}$ for our choice $\sigma_\ell = 1/\log X_\ell$. We define the new event
\begin{equation}\label{eqn:j_Euler_event}
    \mathcal{P}_\ell^{(e)} := \Big\{\sup_{1 \le j \le J}U_j > \frac{T(\ell)}{C\ell^{K(1-2/r)/2}\ell \log \ell}\Big\}
\end{equation}
for some appropriate $C>1$ an absolute constant. Now, this event is independent of $i$ entirely, and we have $\mathbb{P}(\mathcal{Q}_{\ell,1}^{(e)}) \le \mathbb{P}(\mathcal{P}_\ell^{(e)})$. We define $\mathcal{P}_\ell^{(w)}$ similarly using the different weighting. Next, we consider the quantity
\begin{equation*}\begin{split}
    I_{j,\ell}^{(e)} :=& N(\sigma_\ell,y_{j})\int_{\mathbb{R}}\frac{|F_{y_{j}}(1/2+\sigma_\ell+ih)|^2|\zeta(1+2\sigma_\ell+2ih)|^2}{|1/2+\sigma_\ell+ih|^2}dh, \\
    I_{j,\ell}^{(w)} :=&N(\sigma_\ell,y_{j}) \int_\mathbb{R} \frac{|F_{y_{j}}(1/2+\sigma_\ell+ih)|^2}{|\sigma_\ell+ih|^2}dh.
\end{split}
\end{equation*}
In each respective case, we note that
\begin{align*}\Big|\sum_{m \ge 1}\sum_{P(n) \le y}\frac{f(n)}{(m^2n)^s}\Big|^2 =& |\zeta(2s)F_{y}(s)|^2, \qquad \Re(s) > 1/2 \\
\Big|\sum_{P(n) \le y}\frac{f(n)}{n^s}\Big|^2 =& |F_y(s)|^2.
\end{align*}
To allow for the application of Proposition \ref{prop:upper_uniform_REP} in its stated form, we set $s=\sigma+it$ and observe
$$\zeta(2s) \ll \max\Big\{\frac{1}{|2s-1|},\log(2+|t|)\Big\}, \qquad \Re(s) > 1/2.$$
Then define
\begin{equation}\label{eqn:integral_def}
    I_{j,\ell} : = N(\sigma_\ell,y_j)\int_{\mathbb{R}}\frac{|F_{y_j}(1/2+\sigma_\ell+ih)|^2\log^2(2+|h|)}{\sigma_\ell^2+h^2}dh.
\end{equation}
We note $I_{j,\ell}^{(e)} + I_{j,\ell}^{(w)} \ll I_{j,\ell}$, and $I_{j,\ell}$ is in the correct form for Proposition \ref{prop:upper_uniform_REP} to be applied.
\begin{lemma}\label{lem:martingale_euler_products}
    We have $I_{j,\ell}$ is a martingale with respect to the filtration $(\mathcal{F}_{y_{j}})_{j \ge 0}$. 
\end{lemma}
\begin{proof}
    It suffices to check the third condition. Indeed, we have
    \begin{align*}
        \mathbb{E}[I_{j,\ell}|\mathcal{F}_{y_{j-1}}] =& \mathbb{E}\bigg[N(\sigma_\ell,y_{j})\int_{\mathbb{R}}\frac{|F_{y_{j}}(1/2+\sigma_\ell+ih)|^2\log^2(2+|h|)}{\sigma_\ell^2+h^2}dh\Big|\mathcal{F}_{y_{j-1}}\bigg] \\
        =& N(\sigma_\ell,y_{j})\int_{\mathbb{R}}\frac{\mathbb{E}[|F_{y_{j}}(1/2+\sigma_\ell+ih)|^2|\mathcal{F}_{y_{j-1}}]\log^2(2+|h|)}{\sigma_\ell^2+h^2}dh
    \end{align*}
    by Fubini--Tonelli. Applying the conditional expectation, we get
    this is
    \begin{align*}
        = \prod_{y_{j-1} < p \le y_j}(1+p^{-1-2\sigma_\ell})^{-1}&N(\sigma_\ell,y_{j-1})\int_{\mathbb{R}}\mathbb{E}\Big[\prod_{y_{j-1} < p \le y_j}|1+\frac{f(p)}{p^{1/2+\sigma_\ell+ih}}|^2\Big] \\ \times&\frac{|F_{y_{j-1}}(1/2+\sigma_\ell+ih)|^2\log^2(2+|h|)}{\sigma_\ell^2+h^2}dh \\
        =& I_{j-1,\ell}
    \end{align*}
    using with $s = 1/2+\sigma_\ell+ih$
    $$\mathbb{E}\Big[\big|1+f(p)p^{-s}\big|^2\Big] = \frac{1}{2}|1+p^{-s}|^2 + \frac{1}{2}|1-p^{-s}|^2= 1+p^{-1-2\sigma_\ell}.$$
\end{proof}
With this in hand, we can prove the following conditional high moment estimate.
\begin{lemma}\label{lem:conditional_minkowski}
    Let $r \geq 2$ be a fixed integer. Then for $\ell$ sufficiently large and $1 \le j \le J$, we have
    $$\mathbb{E}[(U_{j,\ell}^{(e)})^r|\mathcal{F}_{y_{j-1}}] \ll_r I_{j-1,\ell}^r.$$
    Similarly, we have under the same conditions
    $$\mathbb{E}[(U_{j,\ell}^{(w)})^r|\mathcal{F}_{y_{j-1}}] \ll_r I_{j-1,\ell}^r.$$
\end{lemma}
\begin{proof}
    Our first step in this proof is to apply conditional Minkowski's inequality. Indeed, we have that $(\mathbb{E}[X_j^r|\mathcal{F}_{y_{j-1}}])^{1/r}$ is a norm for $r \geq 1$ for $X_j$ a suitable random variable. Then we apply this to our situation to see
    \begin{equation*}\begin{split}
        (\mathbb{E}[(U_{j,\ell}^{(e)})^r|\mathcal{F}_{y_{j-1}}])^{1/r} :=& \bigg(\mathbb{E}\bigg[\bigg(N(\sigma_\ell,y_j)\int_1^\infty\max_{y_{j-1} \le p < y_j }\Big|\sum_{m \ge 1}\sum_{\substack{n \le z/m^2 \\ P(n) <p}}f(n) \Big|^2\frac{dz}{z^{2+2\sigma_\ell}}\bigg)^r\bigg]\bigg)^{1/r} \\ \le& N(\sigma_\ell,y_j)\int_1^\infty\Big(\mathbb{E}\Big[\max_{y_{j-1} \le p < y_j}\Big|\sum_{m \ge 1}\sum_{\substack{n \le z/m^2 \\ P(n) < p}}f(n)\Big|^{2r} \big|\mathcal{F}_{y_{j-1}}\Big]\Big)^{1/r}\frac{dz}{z^{2+2\sigma_\ell}}.
        \end{split}
    \end{equation*}
    From our observations before, the sum within the integral is a submartingale, so we can apply conditional Doob's $L^r$ inequality Lemma \ref{Lr-Doob} to show this is
    $$\ll_r N(\sigma_\ell,y_j)\int_1^\infty\Big(\mathbb{E}\Big[\Big|\sum_{m \ge 1}\sum_{\substack{n \le z/m^2 \\ P(n) < y_j}}f(n)\Big|^{2r} \big|\mathcal{F}_{y_{j-1}}\Big]\Big)^{1/r}\frac{dz}{z^{2+2\sigma_\ell}}.$$
    At this point, we use Proposition \ref{prop_hyper} to deduce the further upper bound
    $$\ll_r N(\sigma_\ell,y_j)\int_1^\infty\sum_{\substack{k \le z \\\mu^2(k)= 1 \\ p|k \Rightarrow y_{j-1} \le p < y_j}}\tau_{2r-1}(k) \Big|\sum_{m \ge 1}\sum_{\substack{n \le z/km^2 \\ P(n) < y_{j-1}}}f(n)\Big|^{2}\frac{dz}{z^{2+2\sigma_\ell}}.$$
    At this point we use the substitution $v=z/k$, and upper bound the outside sum by
    $$N(\sigma_\ell,y_j)\sum_{p|k \Rightarrow y_{j-1} \le p < y_j}\frac{\mu^2(k)\tau_{2r-1}(k)}{k^{1+2\sigma_\ell}}\int_1^\infty\Big|\sum_{m \ge 1}\sum_{\substack{n \le v/m^2 \\ P(n) < y_{j-1}}}f(n)\Big|^{2} \frac{dv}{v^{2+2\sigma_\ell}}.$$
    The sum over the divisors is then $\prod_{y_{j-1} \le p < y_j}(1+\frac{2r-1}{p}) \ll e^{2r/\ell} \ll_r 1$ as $r$ is fixed. Finally, we apply Proposition \ref{HA1}, so we have in all
    $$(\mathbb{E}[(U_{j,\ell}^{(e)})^r|\mathcal{F}_{y_{j-1}}])^{1/r} \ll_r N(\sigma_\ell,y_j)\int_{-\infty}^\infty \frac{|F_{y_{j-1}}(1/2+\sigma_\ell+ih)|^2\log^2(2+|h|)}{\sigma_\ell^2+h^2}dh \ll_r I_{j-1,\ell}$$
    as required. The case involving the weighted Euler product follows in a near identical manner.
\end{proof}
Before the next lemma, we define the events $\mathcal{I}_{j,\ell} := \{I_{j,\ell} \le \frac{\sqrt{T(\ell)}}{\ell^{K/2}\sqrt{\ell\log\ell}}\}$ and $\mathcal{I}_\ell := \bigcap_{0 \le  j \le J}\mathcal{I}_{j,\ell}$.
\begin{lemma}
    Suppose $r \geq 2$ is an integer. Then $\mathbb{P}(\mathcal{Q}_{\ell,1}^{(e)})$ and $\mathbb{P}(\mathcal{Q}_{\ell,1}^{(w)})$ are summable.  \label{Q*conv}
\end{lemma}
\begin{proof}
By the triangle inequality and the union bound, we have
\begin{align*}
    \mathbb{P}(\mathcal{Q}_{\ell,1}^{(e)}) \le& \mathbb{P}\Big( \mathcal{P}_\ell^{(e)} \bigcap \mathcal{I}_\ell\Big) + \mathbb{P}(\overline{\mathcal{I}_\ell}) \\
    \ll_r& \sum_{1 \le j \le J}\mathbb{P}\Big(\big\{U_{j,\ell}^{(e)} \geqslant \frac{T(\ell)}{\ell^{K(1-2/r)/2}\ell\log \ell} \big\} \bigcap \{\mathcal{I}_{j-1,\ell}\}\Big) + \mathbb{P}(\overline{\mathcal{I}_\ell}). \numberthis \label{Q*ineq}
\end{align*}
We first look at the first sum. Notice that after applying Markov's inequality to the $r$-power, we have
$$\mathbb{P}\Big(\big\{U_{j,\ell}^{(e)} \geqslant \frac{T(\ell)}{\ell^{K(1-2/r)/2}\ell\log \ell} \big\} \bigcap \{\mathcal{I}_{j-1,\ell}\}\Big) \ll_r \frac{\ell^{Kr/2}(\ell\log \ell)^r}{\ell^{K}(T(\ell))^r}\mathbb{E}((U_{j,\ell}^{(e)})^r | \mathcal{I}_{j-1,\ell}).$$
We next use the tower law of conditional expectation noting that $\mathcal{I}_{j-1,\ell} \subseteq \mathcal{F}_{y_{j-1}}$, so we have
$$\mathbb{E}[(U_{j,\ell}^{(e)})^r|\mathcal{I}_{j-1,\ell}] = \mathbb{E}[\mathbb{E}((U_{j,\ell}^{(e)})^r|\mathcal{F}_{y_{j-1}})|\mathcal{I}_{j-1,\ell}] \ll_r \mathbb{E}[I_{j-1,\ell}^r|\mathcal{I}_{j-1,\ell}]$$
using Lemma \ref{lem:conditional_minkowski}. Substituting this back into equation \ref{Q*ineq}, we get
\begin{align*}
    \sum_{1 \le j \le J}\mathbb{P}\Big(\big\{U_{j,\ell}^{(e)} > \frac{T(\ell)}{\ell^{K(1-2/r)/2}\ell\log \ell}\big\}\bigcap\{\mathcal{I}_{j-1,\ell}\}\Big) \ll_r& \sum_{1\le j \le J}\frac{\ell^{Kr/2}(\ell\log \ell)^r\mathbb{E}[I_{j-1,\ell}^r|\mathcal{I}_{j-1,\ell}]}{\ell^{K}(T(\ell))^r} \\
    \ll_r& \sum_{1 \le j \le J} \frac{(\ell\log \ell)^{r/2}}{(T(\ell))^{r/2}\ell^{K}} \ll_r \frac{(\ell\log \ell)^{r/2}}{(T(\ell))^{r/2}},
\end{align*}
which given our choice of $T(\ell)$ is summable in $\ell$.
We now turn our attention to $\mathbb{P}(\overline{\mathcal{I}_{\ell}})$. Then
$$\mathbb{P}(\overline{\mathcal{I}_\ell}) \le \mathbb{P}\Big(\sup_{0 \le j \le J}I_{j,\ell} > \frac{\sqrt{T(\ell)}}{\ell^{K/2}\sqrt{\ell\log \ell}}\Big).$$ We recall by Lemma \ref{lem:martingale_euler_products} that $I_{j,\ell}$ is a martingale in $j$, so $I_{j,\ell}^{1/3}$ is a supermartingale. Consequently, we can apply Doob's maximal inequality, and by Proposition \ref{Lr-Doob}, we have
$$\mathbb{P}\Big(\sup_{0 \le j \le J}I_{j,\ell} > \frac{\sqrt{T(\ell)}}{\ell^{K/2}\sqrt{\ell\log \ell}}\Big) \ll \frac{\ell^{K/6}(\ell\log \ell)^{1/6}\mathbb{E}[I_{0,\ell}^{1/3}]}{(T(\ell))^{1/6}} \ll \frac{(\ell \log \ell)^{1/6}}{(T(\ell))^{1/6}}$$
which converges in $\ell$ given $T(\ell)=\ell^{10}$.
\end{proof}
\subsection{Applying the Caich--Hoeffding inequality}
We now have all of the ingredients to conclude our theorem. We only cover the extended sum since the weighted sum is identical. We have that there exists an absolute constant $C>1$ such that almost surely
$$\frac{V_\ell^{(e)}(x_i)}{x_iC} \ll \ell \log \ell\sup_{\substack{1 \le j \le J \\ 1 \le \frac{\log x_i}{\log y_{j-1}} \le \ell^{2K}}}Q_{j,\ell}^{(e)}(x_i) + \ell^K\sup_{\substack{1 \le j \le J \\ \frac{\log x_i}{\log y_{j-1}} > \ell^{2K}}}Q_{j,\ell}^{(e)}(x_i) + \frac{D_\ell^{(e)}(x_i)}{x_i}.$$
Define the event
\begin{equation}\label{eqn:variance_conditioning}
    \mathcal{V}_{\ell,i}^{(e)} := \Big\{V_{\ell}(x_i;g) \le \frac{Cx_i T(\ell)}{\ell^{K(1-2/r)/2}}\Big\}, \qquad \mathcal{V}_\ell^{(e)}:= \bigcap_{X_{\ell-1} < x_i \le X_\ell}\mathcal{V}_{\ell,i}^{(e)}
\end{equation}
Then we have
$$\mathbb{P}(\overline{\mathcal{V}_\ell^{(e)}}) \le \mathbb{P}(\mathcal{Q}_{\ell,1}^{(e)}) + \mathbb{P}(\mathcal{Q}_{\ell,2}^{(e)}) + \mathbb{P}(\mathcal{D}_\ell^{(e)}),$$
all of the terms of the right hand side of the inequality are summable, so $\mathcal{V}_\ell^{(e)}$ occurs almost surely. Our final step is to apply Proposition \ref{Caich_AH_in}. Define
$$\mathcal{B}_\ell^{(e)} := \Big\{\sup_{X_{\ell-1} < x_i \le X_\ell}|\frac{1}{\sqrt{x_i}}\sum_{\substack{n \le x_i \\ P(\kappa(n)) >y_0}}g(n)| > \ell^{K(1/4+\eta+1/2r)}\Big\}.$$
The same manipulations go through in the case of weighted squarefree Rademacher random multiplicative functions without the $\sqrt{x_i}$ being required in the expressions. We similarly define $\mathcal{B}_\ell^{(w)}$ to $\mathcal{B}_\ell^{(e)}$ with the same changes.
\begin{lemma}\label{lem:gaussian_deviations}
    Both $\mathbb{P}(\mathcal{B}_{\ell}^{(e)})$ and $\mathbb{P}(\mathcal{B}_\ell^{(w)})$ are summable in $\ell$.
\end{lemma}
\begin{proof}
    We appeal to Proposition \ref{Caich_AH_in} to prove this lemma. To this end, we use the event defined in equation (\ref{eqn:variance_conditioning}). By the choice of $T(\ell)$, we know that $\mathbb{P}(\overline{\mathcal{V}_\ell^{(e)}})$ is summable in $\ell$ from previous reasoning. This allows us to write
\begin{align*}
    \mathbb{P}(\mathcal{B}_{\ell}^{(e)}) =& \ \mathbb{P}\Big(\sup_{X_{\ell-1} < x_i \le X_\ell}\frac{|\frac{1}{\sqrt{x_i}}\sum_{\substack{n \le x_i \\ P(\kappa(n)) > y_0}}g(n)|}{\ell^{K(1/4+\eta+1/2r)}} > 1\Big) \\
    \le& \ \mathbb{P}\Big(\bigcup_{X_{\ell-1} < x_i \le X_\ell}\Big(\Big\{\frac{|\frac{1}{\sqrt{x_i}}\sum_{\substack{n \le x_i \\ P(\kappa(n)) > y_0}}g(n)|}{\ell^{K(1/4+\eta+1/2r)}} > 1\Big\} \bigcap \{\mathcal{V}_{\ell}^{(e)}\}\Big)\Big) + \mathbb{P}(\overline{\mathcal{V}_\ell^{(e)}}) \\
    \le& \sum_{X_{\ell-1} < x_i \le X_\ell}\mathbb{P}\Big(\Big\{\frac{|\frac{1}{\sqrt{x_i}}\sum_{\substack{n \le x_i \\ P(\kappa(n)) > y_0}}g(n)|}{\ell^{K(1/4+\eta+1/2r)}} > 1\Big\} \bigcap \{\mathcal{V}_{\ell,i}\}\Big) + \mathbb{P}(\overline{\mathcal{V}_\ell^{(e)}}).
\end{align*}
Applying Proposition \ref{Caich_AH_in}, we have
$$\mathbb{P}\Big(\Big\{\frac{|\frac{1}{\sqrt{x_i}}\sum_{\substack{n \le x_i \\ P(\kappa(n)) > y_0}}g(n)|}{\ell^{K(1/4+\eta+1/2r)}} > 1\Big\} \bigcap \{\mathcal{V}_{\ell,i}^{(e)}\}\Big) \ll \exp\Big(-\frac{C\ell^{K(1+4\eta+2/r)/2}}{T(\ell)\ell^{K(2/r-1)/2}}\Big) \ll \exp\left(-C_1\ell^{K+2K\eta-10}\right)$$
for some $C_1 >0$. Then using $K\eta>24$, we sum over the $x_i$ in the interval $(X_{\ell-1},X_\ell]$. We then have
\begin{align*}
    \mathbb{P}(\mathcal{B}_{\ell}^{(e)} \ll& \sum_{X_{\ell-1} < x_i \le X_\ell} \exp\left(-C_1\ell^{K}\ell^{38}\right) + O(\ell^{-7/6}) 
    \ll \exp(-\ell^K(C_1\ell^{38}-\log 2/\gamma)) + O(\ell^{-7/6}),
\end{align*}
which is summable in $\ell$. Then by the first Borel--Cantelli lemma, $\mathcal{B}_\ell^{(e)}$ occurs at most finitely often almost surely. 
\end{proof}
\printbibliography[title= References]
\end{document}